\documentclass[11pt, oneside]{article}   	
\usepackage{geometry}                		
\usepackage{graphicx}				
\usepackage{amsmath}									
\usepackage{amssymb}
\usepackage{algorithmic}
\usepackage{hyperref}
\newcommand{\pqrs}{\left(
\begin{smallmatrix} 
p & q\\ 
r & s 
\end{smallmatrix}
\right)}

\newcommand{\homographic}[4]{\begin{pmatrix} #1 & #2\\ #3 & #4 \end{pmatrix}}
\newcommand{\bihomographic}[8]{\left(\begin{smallmatrix}#1&#2&#3&#4\\#5&#6&#7&#8\end{smallmatrix}\right)}
\newcommand{\bihom}[8]{{\frac{#1 xy + #2 x + #3 y + #4}{#5 xy + #6 x + #7 y + #8}}}

\newcommand{\abcd}{\left(
\begin{smallmatrix} 
a & b & c & d\\ 
e & f & g & h
\end{smallmatrix}
\right)}

\renewcommand{\:}{\negthickspace:\negthickspace}

\title{Transcendental Functions on Continued Fractions}
\author{Michael J. Collins\\Associate Research Scientist\\Daniel H. Wagner Associates\\Hampton, VA USA\\mjcollins10@gmail.com}
\date{} 

\begin{document}
\maketitle
\begin{abstract}
Gosper  developed an algorithm for performing arithmetic operations on continued fractions (CFs), getting a CF as the result. Straightforward implementation of the algorithm leads to infinite loops on some inputs.
Here we present a modified version of the arithmetic algorithm and prove that it avoids all difficulties with infinite loops.
We then combine CF arithmetic with the ``spigot" algorithm to compute the CF expansions of exponential, logarithmic, and trigonometric functions of CFs.
We have implemented these algorithms in Haskell.
\end{abstract}

\section{Introduction}
A simple continued fraction\footnote{\href{https://en.wikipedia.org/wiki/Simple\_continued\_fraction}{https://en.wikipedia.org/wiki/Simple\_continued\_fraction}} is a (possibly infinite) expression of the form
\[
a_0+\cfrac{1}{a_1+\cfrac{1}{a_2+\cfrac{1}{a_3+\cdots}}}
\]
where the terms $a_i$ are integers. It can be written more compactly as $[a_0, a_1, \cdots]$\footnote{Many authors write this as $[a_0; a_1, \cdots]$}. Here we consider only \emph{regular} simple continued fractions, which means $a_i \geq 1$ when $i>0$. The
properties of continued fractions, and the motivations for studying them, are very well-known \cite{rockett1992continued},
so here we only remind the reader of some notation we will use:
\begin{itemize}
\item A quadratic irrational has a periodic CF expansion. We denote the periodic part by putting it in parentheses,
i.e $[1,(2,3)] = [1,2,3,2,3,2,3,\cdots]$.
\item A rational number has a finite CF expansion $[a_0,a_1, \cdots a_k]$.
It will be convenient to consider a finite CF as ending with an infinite term $a_{k+1}=\infty$, equating $\frac{1}{\infty}$ with zero.
\end{itemize}
 
Gosper (in an appendix to the famous HAKMEM report \cite{hakmem}) 
developed an algorithm for performing arithmetic on two CFs, getting a CF as the result;
the point of course is that we can do this entirely within the CF representation, making no use of floating-point operations. Straightforward implementation of the algorithm leads to infinite loops on some inputs\footnote{\href{https://srossd.com/blog/2020/gosper-1/}{https://srossd.com/blog/2020/gosper-1/}} (section \ref{sec:FailToConverge}).
Here we present a modified version of the arithmetic algorithm, and prove that these modifications avoid infinite loops in the following sense: given any $\varepsilon > 0$, if the correct result is $z=[z_0,z_1,\cdots]$, our algorithm will eventually produce a finite CF $[z_0, \cdots z_{k-1}, z'_k]$ with proof that 
\[
|z - [z_0, \cdots z_{k-1}, z'_k]| < \varepsilon\ .
\]

We then extend the arithmetic algorithm with the techniques of the spigot algorithm for the digits of $\pi$ \cite{Gibbons2016,Rabinowitz2016}, to compute CF representations of $e^x, \log(x)$,  $\cos(x)$, and $\arcsin(x)$, where $x$ is given as a CF. 
Of course, since we can do arithmetic, we could compute any $n^{th}$ degree Taylor approximation of such functions to any degree of accuracy, but we go further;
we develop an iterative algorithm that continually generates terms of the CF, using only as much of the series (and as many terms of $x$) as are required to determine the next output term.

We have implemented these algorithms
in Haskell\footnote{\href{https://github.com/mjcollins10/ContinuedFractions}{https://github.com/mjcollins10/ContinuedFractions}}.

\section{Arithmetic on One CF}
Before we describe how to add or multiply two CFs, we consider the simpler problem of operations combining a single CF with a
rational number $p/q$. It is not in fact necessary to handle this as a special case, but considering this case separately allows us to present the core idea of the algorithm in a straightforward way.

A few examples reveal there is no evident general pattern for transforming the CF terms of irrational $x$ into
the terms of $px/q$ or $x + p/q$, even in the apparently simplest cases:
\begin{eqnarray*}
\sqrt{7} & = & [2, (1, 1, 1, 4)] \\
\sqrt{7}/2 & = & [1, (3, 10, 3, 2)] \\
\sqrt{11} & = & [3, (3, 6)] \\
\sqrt{11}/2 & = & [1, (1, 1, 1, 12, 1, 1, 1, 2)] \\
\pi & = & [3,7,15,1,292,1,1,1, 2, 1, 3,\cdots]\\
\pi + 1/2 & = & [3, 1, 1, 1, 3, 1, 3, 4, 73, 6, 3, 3, 2, 1, 3\cdots]
\end{eqnarray*}
As a motivating example for the general algorithm, we compute
\[[y_0,y_1,y_2 \cdots] = \frac{\pi}{2} = 1.5707963267948966\cdots\ .
\]
Let $\pi = [\pi_0,\pi_1,\cdots]$ (in general we use $x_i$ for the $i^{\mbox {th}}$ CF term of a real number $x$). It will be convenient to have notation for the ``tail" of a continued fraction, so let
\[
r_i = [\pi_{i+1}, \pi_{i+2},\cdots]\ .
\]
Note that all $r_i \geq 1$, and $r_i = \pi_{i+1} + 1/r_{i+1}$. The fact that $\pi_0 = \lfloor \pi \rfloor= 3$
is enough to determine $y_0=\lfloor \pi/2 \rfloor=1$; more precisely
\[
y = \frac{\pi}{2} = \frac{3 + 1/r_0}{2} = \frac{3r_0+1}{2r_0} = 1 + \frac{r_0+1}{2r_0}= 1 + \frac{1}{[y_1,y_2,\cdots]}\ .
\]
Now we have to get $y_1$ from 
\[
[y_1,y_2,\cdots] = \frac{2r_0}{r_0+1}\ .
\]
The mere fact that $r_0>1$ is enough to tell us that $\lfloor \frac{2r_0}{r_0+1} \rfloor = 1$, so we obtain $y_1=1$ and continue
\[
[y_1,y_2,\cdots] = 1 + \frac{r_0-1}{r_0+1} = y_1 + \frac{1}{[y_2,y_3\cdots]}
\]
\[
[y_2,\cdots] = \frac{r_0+1}{r_0-1} \ .
\]
Now the floor of $\frac{r_0+1}{r_0-1}$ ranges from one to infinity as $r_0$ ranges from one to infinity,
so we need to make use of $\pi_1 = 7$; we substitute $r_0 = 7+1/r_1$ to get
\[
[y_2,\cdots] = \frac{8 + 1/r_1}{6 + 1/r_1} = \frac{8r_1 + 1}{6r_1 + 1} = 1 + \frac{2r_1}{6r_1 + 1}= 1 + \frac{1}{[y_3,\cdots]}\ .
\]
Thus $y_2=1$.  We also get $y_3 = 3$ from
\[
[y_3,\cdots] =  \frac{6r_1 + 1}{2r_1} = 3 + \frac{1}{2r_1} = 3 + \frac{1}{[y_4,\cdots]}\ .
\]

Going further will require substituting $r_1 = 15+1/r_2$ into $[y_4,\cdots] = 2r_1$.
It is now clear that we will generate terms of $y$ by repeatedly determining the integer part (i.e. floor) of expressions of the form
\begin{equation}\label{eq:homographic}
\frac{px+q}{rx+s}
\end{equation}
where $x > 1$, the coefficients are integral, and the continued fraction expansion of $x$ is known. Functions of this form are called \emph{homographic}.
We will identify a homographic function with the matrix
$M=\left(
\begin{smallmatrix} 
p & q\\ 
r & s 
\end{smallmatrix}
\right)$
and may write (\ref{eq:homographic}) as $M(x)$. We can assume without loss of generality that the determinant $ad-bc$ is nonzero, and $\gcd(a,b,c,d)=1$.

If $-s/r \leq 1$ then the denominator cannot be zero, and $M(x)$ is between the min and max of $\{p/r, (p+q)/(r+s)\}$; if
\[
\left\lfloor \frac{p}{r}\right\rfloor = \left\lfloor \frac{p+q}{r+s}\right\rfloor
\]
then this common value is $\lfloor M(x) \rfloor$; we say that $M$ is \emph{unambiguous}.

We always have
\begin{equation}\label{eq:invar}
M([x_i,x_{i+1},\cdots]) = [y_j,y_{j+1},\cdots]
\end{equation}
where $[x_i,x_{i+1},\cdots]$ is the unread part of the input and $[y_0,\cdots y_{j-1}]$ is the output prefix already computed. 
We define two transformations on $M$ which preserve (\ref{eq:invar}). When $M$ is ambiguous, we do not yet have enough information to determine $y_j$, and must $\emph{ingest}$ the next term of $x$. If this term is $x_i=k$ we make the substitution $x \leftarrow k + 1/x$, leading to the definition
\[
\mbox{ingest}(k,\pqrs) = \homographic{q+kp}{p}{s+kr}{r}\ .
\]
When $M$ is unambiguous with $k=\left\lfloor \frac{p}{r}\right\rfloor$, we can \emph{produce}
the next output term $y_j = k$; the rest of the output, $[y_{j+1},y_{j+2},\cdots]$, is the continued fraction expansion of
\[
\frac{1}{\frac{px+q}{rx+s}-k} = \frac{rx+s}{(p-kr)x + (q-ks)}
\] 
so we define 
\[
\mbox{produce}(k,\pqrs) = \homographic{r}{s}{p-kr}{q-ks}\ .
\]
Note that these two transformations do not change the absolute value of the determinant of $M$.

If we started with rational $x$, we will eventually reach $x_i = \infty$; ingesting infinity returns the limit ${\left(
\begin{smallmatrix} 
0 & p\\ 
0 & r
\end{smallmatrix}
\right)}$.
This is the constant function at the rational number $p/r$, so there will never be any need to ingest the (nonexistent) subsequent 
terms of $x$. If we proceed as in the examples above, we will compute the CF representation of $p/r$ as the last part of the finite
expansion of $y$, ending with a production that subtracts the final term $y_k$ from the expression $\frac{y_k}{1}$; this leads to
$M=\left(\begin{smallmatrix}0 & 1 \\ 0 & 0\end{smallmatrix}\right)$, which we may read as infinity.

\begin{figure}\label{fig:oneCFarith}
\begin{algorithmic}
\STATE{}\COMMENT{input $x = [x_0,x_1,\cdots]$}
\STATE{}\COMMENT{input $M = \pqrs$}
\STATE{}\COMMENT{output $y$: CF expansion of $M(x)$}
\STATE{$i \gets 0$} \COMMENT{index of the next term of $x$ we will read}
\STATE{$j \gets 0$} \COMMENT{index of the next term of $y$ we will generate}
\STATE{$M \gets \mbox{ingest}(x_0, M)$} \COMMENT{updating $M$ updates variables $p,q,r,s$}
\STATE{$i \gets 1$}
\WHILE {$M \neq \left(\begin{smallmatrix}0 & 1 \\ 0 & 0\end{smallmatrix}\right)$} 
	\WHILE{$M$ is ambiguous}
                   \STATE{$M \gets \mbox{ingest}(x_i, M)$}
                   \STATE{$i \gets i+1$}
           \ENDWHILE
           \STATE{$y_j \gets \left\lfloor \frac{p}{r}\right\rfloor$}
          \STATE{$M \gets \mbox{produce}(y_j,M)$}
          \STATE{$j \gets j+1$}
\ENDWHILE
\end{algorithmic}
\caption{Algorithm for arithmetic on a CF and a rational number}
\end{figure}

So now we have an algorithm for producing the CF expansion $[y_0,y_1,\cdots]$ of $\frac{px+q}{rx+s}$,
assuming the expansion of $x$ is known.
To compute $x/k$ we would start with $M=\left(
\begin{smallmatrix} 
1 & 0\\ 
0 & k 
\end{smallmatrix}
\right)$; 
to get $kx$ we would start with $M=\left(
\begin{smallmatrix} 
k & 0\\ 
0 & 1 
\end{smallmatrix}
\right)$; while $x+\frac{j}{k}$ would be $M=\left(
\begin{smallmatrix} 
k & j\\ 
0 & k 
\end{smallmatrix}
\right)$. Note that computing any homomorphic expression $\pqrs$ is no more difficult than a single arithmetic operation.

\subsection{Proof of Termination}\label{sec:oneCF}
The outer loop ``while $M \neq \infty$" is an infinite loop if $\mathbf{x} = [\mathbf{x}_0, \mathbf{x}_1,\cdots]$ is irrational
(here we use boldface $\mathbf{x}$ for the real-valued input, as distinct from the symbolic variable $x$ in the homographic
expression). This loop can be implemented quite directly in a language like Haskell\footnote{\href{https://haskell.org}{haskell.org}} with
lazy evaluation\cite{hutton2007programming}; such languages support conceptually infinite data structures, which are a perfect fit for working with continued
fractions. Under lazy evaluation, defining a list-type variable $y$ to be the result of an algorithm just associates $y$ with the finite recursive
expression that defines the algorithm; the number of elements of $y$ actually generated will be only what is needed by subsequent computations.

The inner loop ``while $M$ is ambiguous" \emph{is} guaranteed to terminate even on irrational $\mathbf{x}$. As
we iterate through the inner loop, we are repeatedly rewriting the initial $M = \pqrs$
with $[\mathbf{x}_0, \cdots \mathbf{x}_i, x]$ in place of $x$. The upper and lower bounds on
$[\mathbf{x}_0, \cdots \mathbf{x}_i, x]$ approach $\mathbf{x}$, so the upper and lower bounds of the matrix approach
$\frac{p\mathbf{x}+q}{r\mathbf{x}+s}$.
Since this limit is irrational\footnote{because we have ruled out the degenerate case of zero determinant}, eventually the bounding interval contains no integers; at this point $M$  will be unambiguous
 and we will produce a term.
 
 Liardet and Stambul \cite{LSAlgebraic1998} go beyond this basic observation to prove a bound (dependent on the determinant of $M$) on the maximum number of iterations needed to get the next term. Miska, Murru, and Romeo \cite{miska2025} give performance results suggesting that, almost always, four input terms are sufficient to get the next output.

\subsection{Visualizing the Algorithm}
Since $\mbox{ingest}(k,\pqrs)$ returns a matrix whose right column is
$\left( \begin{smallmatrix} p \\  r \end{smallmatrix}\right)$,
and $\mbox{produce}(k,\pqrs)$ returns a matrix whose top row is $(r\ s)$, we can visualize the progress of the algorithm
as moving a $2 \times 2$ window through a two-dimensional grid of integers.
We start with $M$ in the upper right corner; ingest moves left and produce moves down.
We put the $x_i$ on the top row, at the positions where we ingest, and $y_j$ in the right column, at the positions where we produce.
Here is such a representation of our earlier example of $\pi/2$, taken further to use $\pi = [3,7,15, 1\cdots]$ to obtain $\pi/2 = [1,1,1,3,31,\cdots]$.

\[\label{eq:2dgrid}
\begin{matrix}
 & & & 1& 15&       7&  3& &\\
 & & &  &        &     3&   1& 0&\\
 & & &  &        &     2&    0& 2& 1\\
 & & &  &        8&   1&    1&   &  1\\
 & & &  &        6&   1&    -1&  & 1\\
 & & 32& 30&  2& 0&        &  &  3\\
 & &  1&   1&  0&  1&        &  & 31\\
 & &   1&  -1&   &    &        &  &
\end{matrix}
\]

\section{Arithmetic on Two Continued Fractions}
We now turn to adding or multiplying two CFs. Conceptually, this is hardly any different from what we have already done;
but it will raise some implementation issues that require care. As a motivating example, let
\[
z = \pi + \sqrt{2} = 4.555806215962888\cdots
\]
 with $[y_0, y_1,\cdots] = \sqrt{2} = [1,(2)]$. We can write this sum as
\[
3 + \frac{1}{[\pi_1,\pi_2\cdots]} + 1 + \frac{1}{[y_1,\cdots]} = \frac{ 4[\pi_1,\cdots][y_1,\cdots] + [\pi_1,\cdots]  + [y_1\cdots]}{ [\pi_1,\cdots][y_1\cdots] }\ .
\]
Substituting $[\pi_1,\cdots] = 7 + 1/[\pi_2,\cdots]$ and $[y_1,\cdots] = 2 + 1/[y_2,\cdots]$ leads (after a lot of high-school algebra) to 
\[
z = \frac{65[\pi_2,\cdots] [y_2,\cdots] + 29[\pi_2,\cdots] +9[y_2,\cdots] +4}{14[\pi_2,\cdots] [y_2,\cdots] + 7[\pi_2,\cdots] +2[y_2,\cdots] +1}
\]
The integer parts of $\frac{65}{14}, \frac{29}{7}, \frac{9}{2}$ and  $\frac{4}{1}$ are all 4; therefore $z_0 = \lfloor z \rfloor =4$.
This conclusion follows easily from the mediant inequality: if $\frac{a}{c} < \frac{b}{d}$ with all terms positive, then
\[
\frac{a}{c} < \frac{a+b}{c+d}< \frac{b}{d}
\] 
The next term $z_1$ will be the floor of
\[
\frac{1}{z-4} =  \frac{14[\pi_2,\cdots] [y_2,\cdots] + 7[\pi_2,\cdots] +2[y_2,\cdots] +1}{9[\pi_2,\cdots] [y_2,\cdots] + [\pi_2,\cdots] +[y_2,\cdots]}
\]
into which we can substitute  $[\pi_2,\cdots]=15+\frac{1}{[\pi_3,\cdots]},  [y_2,\cdots]=2+\frac{1}{[y_3,\cdots]}$ and so on.

Similarly, the product $\pi\sqrt{2}=4.442882938158366\cdots$ is
\[
(3 + \frac{1}{[\pi_1,\cdots]})(1 + \frac{1}{[y_1\cdots]}) 
  = \frac{ 3[\pi_1,\cdots][y_1\cdots] + 3[\pi_1,\cdots] + [y_1\cdots] +1 }{[\pi_1,\cdots][y_1\cdots] }
\]
into which we can make the same substitutions. 

In general, computing terms of a sum or product of CFs will require finding the floors of two-variable expressions of the form
\begin{equation}\label{eq:bihomographic}
\frac{axy + bx + cy + d}{exy + fx + gy + h}
\end{equation}
where $x$ and $y$ independently vary from 1 to $\infty$. Such an expression is called \emph{bihomographic}, and will be
represented by the matrix $\abcd$. To determine bounds on the floor of (\ref{eq:bihomographic}), it is convenient to make the
substitutions $\hat{x} = x-1, \hat{y} = y-1$, and consider
\[
\frac{a\hat{x}\hat{y} + (a+b)\hat{x} + (a+c)\hat{y} + (a+b+c+d)}{e\hat{x}\hat{y} + (e+f)\hat{x} + (e+g)\hat{y} + (e+f+g+h)}
\]
as $\hat{x}, \hat{y}$ range independently from \emph{zero} to infinity. If the denominator cannot be zero, the floor is always
between the minimum and maximum of
\[
\left\{ \frac{a}{e},\frac{a+b}{e+f},\frac{a+c}{e+g},\frac{a+b+c+d}{e+f+g+h} \right\}
\]
where we may ignore fractions with numerator and denominator both zero.

When the floor (i.e. the next term of output) is known, we produce output $z_j=k$ and the next bihomographic expression is
\[
\left(\frac{axy + bx + cy + d}{exy + fx + gy + h} - k\right)^{-1} = \frac{exy + fx + gy + h}{(a-ke)xy + (b-kf)x + (c-kg)y + (d-kh)}
\]
so we define the corresponding operation on matrices
\[
\mbox{produce}(k, \abcd) = \bihomographic{e}{f}{g}{h}{(a-ke)}{(b-kf)}{(c-kg)}{(d-kh)}
\]
 If the floor is not determined (in particular when the denominator might be zero), we must narrow the range by ingesting the next term of either $x$ or $y$.
 If we use $x = [s,x_{k+1},x_{k+2}\cdots]$, we make the substitution $x \leftarrow s + 1/x$ to get
\[
\frac{(sa+c)xy + (sb+d)x + ay + b}{(ae+g)xy + (af+h)x + ey + f}
\]
so we define
\[
\mbox{ingest\_x}(s, \abcd) = \bihomographic{(sa+c)}{(sb+d)}{a}{b} {(ae+g)}{(af+h)}{e}{f}
\]
with the analogous 
\[
\mbox{ingest\_y}(s, \abcd) = \bihomographic{(sa+b)}{a}{(sc+d)}{c}{(se+f)}{e}{(sg+h)}{g}
\]
for the result of substituting $y \leftarrow s + 1/y$.

If one or both inputs are rational, we will eventually ingest infinity, which we define by taking the limit. In most cases
\[
\mbox{ingest\_x}(\infty, \abcd) = \bihomographic{0}{0}{a}{b} {0}{0}{e}{f}
\]
\[
\mbox{ingest\_y}(\infty, \abcd) = \bihomographic{0}{a}{0}{c}{0}{e}{0}{g} \ .
\]
However, this definition of $\mbox{ingest\_x}$ is valid only when the numerator and denominator both contain $x$ terms
(i.e. only if $a,b$ are not both zero and $e,f$ are not both zero), with an analogous requirement for $\mbox{ingest\_y}$.
If the matrix is of the form $\frac{bx+d}{fy+h}$ then the limit as $x \to \infty$ is infinity\footnote{with the caveat that the denominator cannot be zero, which occurs if we divide by zero; see section \ref{sec:anomalous}.}
and the limit as $y \to \infty$  is zero.

\begin{figure}\label{fig:twoCFarithPrelim}
\begin{algorithmic}
\STATE{}\COMMENT{input $x,y = [x_0,x_1,\cdots], [y_0,y_1,\cdots]$}
\STATE{$i \gets 0$} \COMMENT{index of the next term of $x$ and $y$ we will read}
\STATE{$j \gets 0$} \COMMENT{index of the next term of $z$ we will generate}
\STATE{$M \gets \abcd$} \COMMENT{the initial matrix, i.e. the expression whose CF we will compute}
\STATE{$M \gets \mbox{ingest\_y}(y_0, \mbox{ingest\_x}(x_0, M))$}\COMMENT{updating $M$ updates variables $a,\cdots, h$}
\STATE{$i \gets 1$}
\WHILE {$M \neq \infty$} 
	\WHILE{$M$ is ambiguous}
                   \STATE{$M \gets \mbox{ingest\_y}(y_i, \mbox{ingest\_x}(x_i, M))$}
                   \STATE{$i \gets i+1$}
           \ENDWHILE
           \STATE{$z_j \gets \lfloor M(x,y)\rfloor$}
          \STATE{$M \gets \mbox{produce}(z_j,M)$}
          \STATE{$j \gets j+1$}
\ENDWHILE
\end{algorithmic}
\caption{Preliminary algorithm for arithmetic on two CFs}
\end{figure}

We now have a preliminary algorithm for arithmetic, which is essentially the original HAKMEM algorithm in more explicit notation. 
This algorithm will in fact require substantial modification, but we include it as a summary
of the key ideas. It is not much different from the single-CF case. Computations of $x+y, x-y,xy$, and $x/y$ begin with matrices
$\bihomographic{0}{1}{1}{0}{0}{0}{0}{1}, \bihomographic{0}{1}{-1}{0}{0}{0}{0}{1},
\bihomographic{1}{0}{0}{0}{0}{0}{0}{1}$ and $\bihomographic{0}{1}{0}{0}{0}{0}{1}{0}$.

Note that it is not at all necessary to ingest terms of both inputs in lockstep as we have here;
Gosper suggests heuristics for accelerating the algorithm by choosing the input term most likely to get us to the next production step
more quickly. This appears to be a worthwhile subject for further work.

It is interesting that the algorithms for all arithmetic operations are identical, except for the initial $M$. The reason, one might say, is
that division is the hardest arithmetic operation, and the definition of continued fractions has division built-in to everything.

\subsection{Failure to Converge}\label{sec:FailToConverge}
Our preliminary algorithm demonstrates the basic idea of CF arithmetic, but is insufficient because it can fail to converge;
we might endlessly ingest terms of both $x$ and $y$ without ever obtaining a bihomographic expression whose floor is known. The
simplest case is multiplying $\sqrt{2} = [1,(2)]$ by itself. The first iteration, ingesting 1 from $x$ and from $y$, gives
\[
z = \sqrt{2}*\sqrt{2} = [z_0,z_1, \cdots] = \frac{xy+x+y+1}{xy} \ . 
\]
Repeatedly ingesting 2 from both $x$ and $y$ yields
\begin{equation}\label{eq:bound2}
[z_0,z_1, \cdots]  =  \bihom{9}{3}{3}{1}{4}{2}{2}{1}
\end{equation}
\begin{equation}\label{eq:bound3}
[z_0,z_1, \cdots]  =  \bihom{49}{21}{21}{9}{25}{20}{10}{4}
\end{equation}
\begin{equation}\label{eq:bound4}
[z_0,z_1, \cdots]  =  \bihom{289}{119}{119}{49}{144}{60}{60}{25}
\end{equation}
and so on. The matrix (\ref{eq:bound2}) takes values ranging from $\frac{16}{9}$ to $\frac{9}{4}$ (i.e. 1.777 to 2.25)
as $x$ and $y$ range from one to infinity;
similarly matrix (\ref{eq:bound3}) is between $\frac{49}{25}$ and $\frac{100}{49}$ (i.e. 1.96 to 2.0408).
And the last ranges from $\frac{576}{289}$ to $\frac{289}{144}$, which is about 1.993 to 2.00694.

The source of the problem is now clear: we should have $z=[z_0,z_1]=[2,\infty]$, but no finite number of terms of $x$ and $y$ can ever suffice to determine that
$z_0=2$. It is always possible that one input is in fact less than $\sqrt{2}$, so it is always possible that $z<2$ and $z_0=1$.

Our approach to the problem (see \cite{Lester01} for an entirely different approach)
applies basic notions from interval arithmetic (\cite{doi:10.1137/1.9780898717716,Mayer+2017}).
We separate the internal representation of a CF as a data structure from the mathematical notion of a  sequence of integer terms. The
internal representation will contain explicit terms (i.e. $z_k=n$) whenever such can be determined, but might also contain arbitrary
bounds on the tail $[z_k, z_{k+1}\cdots]$, such as the  sequence of bounds we derived for $\sqrt{2}*\sqrt{2}$. This is a direct
generalization, because the explicit term  $z_k=n$ is equivalent to the bound 
\[
n \leq [z_k, z_{k+1}\cdots] < n+1 \ .
\]
Therefore, when a CF is constructed explicitly from a known sequence of integer terms $[a_0, a_1, \cdots]$,
it will be represented internally as the sequence of half-open intervals
\[
[[a_0, a_0+1),  [a_1, a_1+1), \cdots]\ .
\] 
At the opposite extreme, the internal representation for $(\sqrt{2})^2$ will consist entirely of ever-tighter bounds on the first term, i.e. the sequence of nested intervals
\[
\left[ \left[\frac{16}{9},\frac{9}{4}\right), \left[\frac{49}{25},\frac{100}{49}\right), \left[\frac{576}{289},\frac{289}{144}\right), \cdots\right ]\ .
\] 

The general idea (elaborated more fully in section \ref{sec:finalAlgorithm}) is that we will generate output on \emph{every}
 iteration of the arithmetic algorithm; when we cannot produce the next term of $z$, we output upper and lower bounds on the current
 homographic matrix, leaving the matrix unchanged. Reading such output, there can be no ambiguity about the term to which a bound 
 applies; the first bound applies to $[z_0,z_1\cdots]$, and as soon as we encounter an explicit term $[a,a+1)$ we know that $z_k=a$ and the next bound will apply to $[z_{k+1}, z_{k+2}\cdots]$. 
 This implies that any finite prefix of the output sequence will determine a finite prefix $[z_0 \cdots z_{k-1}]$ of $z$, followed by a bound
 $[\ell_k,u_k)$ on the tail $[z_k, z_{k+1}\cdots]$. The exact value of $z$ is between $[z_0 \cdots z_{k-1}, \ell_k]$ and
 $[z_0 \cdots z_{k-1}, u_k]$, if we extend the usual CF notation by allowing a rational number, instead of an integer, as the last term.

Storage of a long sequence of intervals is not an issue for implementation;
successive bounds on the same term are always tighter, so there is no need to retain earlier bounds when a new one is obtained.

A sequence of intervals will represent a real number if every ambiguous interval $I$ is immediately followed
either by another ambiguous interval $I' \subset I$, or by an integer term $a$ with $a \in I$;
if the sequence ends with an infinite sequence of ambiguous intervals, the intersection of that infinite sequence must be a singleton $\{a\}$.
We call such a sequence \emph{valid}.
It will be easy to verify that, if $x$ and $y$ are valid, the result of the arithmetic algorithm $z = M(x,y)$ is valid.
The question of validity will become somewhat more subtle in section \ref{sec:transcendental} when we consider transcendental functions of a CF.

\subsection{Extracting Terms}
Note that our interval representation of CF arithmetic makes it impossible to simply ask a question like ``what are the first
five terms of $z = xy$?" An explicit list of terms can only be relative to a required degree of accuracy $\varepsilon$. Our arithmetic
algorithms consume and return theoretically infinite sequences of intervals; extracting a finite sequence of terms is a separate
``post-processing" step, carried out only when a numerical approximation is required.
To obtain such an approximation, we find a prefix of
the output interval sequence which is long enough to give
\[
 \bigl| [z_0, \cdots z_{k-1}, \ell_k] - [z_0, \cdots z_{k-1}, u_k] \bigr| \leq \varepsilon \ .
\]
Then our rational approximation is $[z_0, \cdots z_{k-1}, z'_k]$, where $z'_k$ is some integer in the interval
$[\ell_k,u_k)$. We might as well always take $z'_k = \lfloor u_k \rfloor$, which is the only integer in the interval if the interval is small. 

In the case of $z=(\sqrt{2})^2$, any $\varepsilon$ will yield the same result, $z = [2]$. More generally, whenever $M(\mathbf{x},\mathbf{y}) = [z_0, \cdots z_k]$ is rational, the rational approximation will eventually be $[z_0, \cdots z_k]$. As in section (\ref{sec:oneCF}), this follows simply from the fact that upper and lower bounds on (\ref{eq:bihomographic}) will approach $z$ as bounds on $x$ and $y$ approach $\mathbf{x}$ and $\mathbf{y}$.

\subsection{An Algorithm for CF Arithmetic}\label{sec:finalAlgorithm}
We now present the full details of the algorithm suggested in the previous section.  We augment the matrix $M$ with a pair of intervals
$\{I_x, I_y\} = \{[x_{\ell},x_{u}), [y_{\ell}, y_{u})\}$, the current bounds on the remaining tails of the inputs $x,y$. These
intervals start at $(-\infty, \infty)$, since the first term could be any integer.

When we read an interval of the form $[a,a+1)$ from
$x$ or from $y$, we modify $M$ as before with $M \leftarrow \mbox{ingest\_x}(a,M)$ or $M \leftarrow \mbox{ingest\_y}(a,M)$,
and we change that variable's bound to $[1, \infty)$. When we read an \emph{ambiguous} interval $[\ell, u)$ (i.e. a bound that
does not determine the next term of the input), we replace the bound $I_x$ or $I_y$  with $[\ell, u)$, leaving $M$ unchanged.

When $M$ is unambiguous with $k = \lfloor M(x,y) \rfloor$, we output $[k,k+1)$ and modify $M$ as before with $M \leftarrow \mbox{produce}(k,M)$;
the intervals $I_x$ and $I_y$ do not change. Otherwise we leave $M$ unchanged and output the interval from min to max of $M(x,y)$ subject to $x \in I_x, y \in I_y$:
denote this interval $\rho\left(M, I_x, I_y\right)$.
These bounds are easily computed, see section (\ref{sec:ComputeBounds}).
To retain unlimited accuracy, the bounds must be computed with exact rational arithmetic.

\subsubsection{Recursive Definition}
We can express the algorithm as a recursive function,
which is essentially the same as our Haskell implementation.
We first define versions of $\mbox{ingest\_x},\mbox{ingest\_y}$ that consider the bounds on $x$ and $y$. Given an explicit term $[s,s+1)$,
\[
\mbox{ingBd\_x}([s,s+1), M, I_x, I_y) = 
(\mbox{ingest\_x}(s,M), (1,\infty),I_y)
\]
while for an ambiguous interval $I_0$,
\[
\mbox{ingBd\_x}(I_0, M, I_x, I_y) = 
(M, I_0,I_y)
\]
with a similar definition for \mbox{ingBd\_y}.

We then define
$\mbox{arith}(x,y, M, I_x, I_y)$
to be the representation of $M(x,y)$ as a list of intervals, where $x,y$ are lists of intervals and $M$ is bihomographic. We use the standard recursive notation $x = x_0\:x_s, y = y_0\:y_s$ to denote the head and tail of a list\footnote{i.e. $a\:[] = [a]$ and $a\:[b,c,\cdots] = [a,b,c,\cdots]$}.  If $\rho(M,I_x,I_y)$ determines $a$ as the next output term, then
\[
\mbox{arith}(x,y, M, I_x, I_y) = (a,a+1)\:\mbox{arith}(x,y,  \mbox{prod}(a,M) ,  I_x,  I_y)
\]
while if $\rho(M,I_x,I_y)$ is ambiguous
\[
\mbox{arith}(x,y, M, I_x, I_y) = \rho(M,I_x,I_y)\:\mbox{arith}(x_s,y_s,  \mbox{ingBd\_y}( y_0, (\mbox{ingBd\_x}(x_0,M,I_x,I_y)) ))
\]
If $M=\infty$ (i.e. the denominator is zero), we have exhausted a pair of rational inputs to produce a rational result, and can terminate the recursion with $\mbox{arith}(x,y, M, I_x, I_y) = [\infty]$.

\begin{figure}
\begin{algorithmic}\label{fig:twoCFarith}
\STATE{}\COMMENT{input $x,y$: CFs represented as lists of intervals}
\STATE{}\COMMENT{input $M$: bihomographic matrix}
\STATE{}\COMMENT{output $z$: CF expansion of $M(x,y)$, represented as a list of intervals}
\STATE{$i \gets 0$} \COMMENT{index of the next element of $x$ and $y$ we will read}
\STATE{$j \gets 0$} \COMMENT{index of the next element of $z$ we will generate}
\STATE{$I_x,I_y \gets (-\infty, \infty),(-\infty, \infty)$} \COMMENT{current bounds on remaining parts of $x,y$}
\WHILE {$M \neq \infty$}
    \STATE{}\COMMENT{read input}
     \IF{$x_i = [a,a+1)$ (where $a$ is an integer)} 
         \STATE{$M \gets \mbox{ingest\_x}( a, M)$}
         \STATE{$I_x \gets [1,\infty)$}
     \ELSE    
          \STATE{$I_x \gets x_i$}
     \ENDIF
          \IF{$y_i = [a,a+1)$ (where $a$ is an integer)} 
         \STATE{$M \gets \mbox{ingest\_y}( a, M)$}
         \STATE{$I_y \gets [1,\infty)$}
     \ELSE    
          \STATE{$I_y \gets y_i$}
     \ENDIF
     \STATE{$i \gets i+1$}
     
     \STATE{}\COMMENT{produce output}
     \IF{$(M,I_x,I_y)$ is unambiguous with $\lfloor M(x,y) \rfloor = a$ for all $x\in I_x, y\in I_y$}
     \STATE{$M \gets \mbox{produce}(a, M)$}
          \STATE{$z_j \gets [a,a+1)$}    
     \ELSE
         \STATE{$z_j \gets \rho(M,I_x, I_y)$}
     \ENDIF
     \STATE{$j \gets j+1$}
\ENDWHILE
\end{algorithmic}
\caption{Algorithm for arithmetic on two CFs (without optimizations)}
\end{figure}

\subsection{Computing the range of $M$}\label{sec:ComputeBounds}
To compute $\rho\left(\abcd, [x_{\ell}, x_u),  [y_{\ell}, y_u)\right)$:
\begin{itemize}
\item Let $\delta_x, \delta_y = x_u - x_{\ell}, y_u - y_{\ell}$
\item Substitute $x \leftarrow x_{\ell} + \frac{\delta}{x'+1}$ and $y \leftarrow y_{\ell} + \frac{\delta}{y'+1}$.
Now $x$ goes from $x_{\ell}$ to $x_u$ as $x'$ goes from 0 to infinity, similarly for $y,y'$
\item Let the coefficients of the resulting matrix be $\bihomographic{a'}{b'}{c'}{d'}{e'}{f'}{g'}{h'}$
\end{itemize}
The range of $M$ subject to the bounds will be the min and max of 
\[
\left\{ \frac{a'}{e'},\frac{b'}{f'},\frac{c'}{g'},\frac{d'}{h'} \right\}
\]
assuming no sign change in the denominator as $x'$ and $y'$ go from 0 to infinity; fractions with numerator and denominator both $0$ can be ignored.
If the denominator can be zero, $\rho\left(\abcd, [x_{\ell}, x_u),  [y_{\ell}, y_u)\right) = [-\infty,\infty)$

\subsection{Possible Optimization}
When we have to read the next interval from $x$ or $y$, and we see an ambiguous interval at the head of the list, it is not necessary
to use that interval; we may skip it and go on to the next element of the list, which will give a tighter bound on the next term.
Any such heuristic must guarantee that, given an infinite list of ambiguous input intervals, we eventually read something
and recompute the range of $M$. Similarly, the
algorithm does not need to output every ambiguous interval of $z$; but  we must guarantee that we eventually output something,
even if we are generating an infinite list of ambiguous intervals. The most obvious idea is to read and generate ambiguous intervals only when the width of the interval is smaller than some threshold.

\subsection{Division by Zero}\label{sec:anomalous}
The limit of $\frac{axy + bx + cy + d}{exy + fx + gy + h}$ as $x$ or $y$ approaches infinity is not necessarily a homographic matrix; for example
\[
\lim_{y \to \infty} \frac{y}{x-3} = \begin{cases}
                                                          +\infty & \mbox{if } x > 3 \\
                                                          \mbox{undefined} & \mbox{if } x=3 \\
                                                          -\infty & \mbox{if } x < 3
                                                         \end{cases}
\]
But ingesting infinity can only mean that we have reached the end of a finite rational input $\mathbf{y}$. At a previous step we
must have ingested an explicit integer term $\alpha$. The previous matrix, which generated $\frac{y}{x-3}$ by substituting $y \leftarrow \alpha + 1/y$, can be found by applying the inverse substitution $y \leftarrow 1/(y-\alpha)$ to get
\[
\frac{1}{(y-\alpha)(x-3)} = \frac{1}{xy - \alpha x -3y +3\alpha}
\]
which means we are trying to divide by zero; we are substituting $y = \alpha + \frac{1}{\infty} = \alpha$. The invalid computation is revealed by the anomalous limit.
Such anomalous limits can only arise from division by zero.

\section{Transcendental Functions}\label{sec:transcendental}
\subsection{Exponentials} We now consider how to compute the CF terms of 
\[
e^x = \sum_{k=0}^{\infty} \frac{x^k}{k!}
\]
where $x$ is given as a CF, possibly with ambiguous intervals. Let
\[
M_{n,x}(y) = 1 + \frac{xy}{n} \ .
\]
This is a bihomographic expression with matrix $\bihomographic{1}{0}{0}{n}{0}{0}{0}{n}$. The $n^{th}$ degree Taylor approximation,
$1 + x + x^2/2 +\cdots + x^n/{n!}$, can be written using Horner's rule as
\[
 1 + x \left(1 + \frac{x}{2}\left( 1 + \frac{x}{3}\left(1 + \cdots \left(1 + \frac{x}{n}\right)\right)\right) \cdots \right) = M_{(1,x)}( M_{(2,x)}( \cdots M_{(n,x)}(1) \cdots ) )
 \]
 so in the limit
 \begin{equation}\label{eq:spigot}
 e^x = 1 + x \left(1 + \frac{x}{2}\left( 1 + \frac{x}{3}\left(1 + \cdots \right. \right.  \right.   \ .
 \end{equation} 
 Now define
 \begin{equation}\label{eq:yn}
 y_n(x) = 1 + \frac{x}{n} + \frac{x^2}{n(n+1)} + \cdots = \lim_{m\to\infty} M_{(n,x)}( M_{(n+1,x)}( \cdots M_{(m,x)}(1) \cdots ) ) \ .
  \end{equation}
Thus $y_n(x) = M_{(n,x)}(y_{n+1}(x))$, and $y_1(x) = e^x$.
Assume that $x \geq 0$ (otherwise we can use the identity $e^{-x} = 1/e^x$); we have
 \[
 y_n(x) - 1 \leq \frac{e^x - 1}{n} \ .
 \]
  We can always obtain arbitrarily tight rational bounds $x_{\ell} \leq x \leq x_{h}$ from a finite prefix of the CF representation of $x$; since $e < 3$  we have the rational bounds
 \begin{equation}\label{eq:ynInterval}
  1+\frac{x_{\ell}}{n} \leq y_n(x) \leq  1 + \frac{3^{\lceil x_{h}\rceil}-1}{n} \ .
 \end{equation}
 If this interval is ambiguous, it can be the first interval of the CF representation of $y_n(x)$.
The subsequent intervals are the CF representation of $M_{(n,x)}(y_{n+1}(x))$,
which is a bihomographic expression that can be evaluated with the arithmetic algorithm;
the algorithm will read
\[
\left(1+\frac{x_{\ell}}{n+1} , 1 + \frac{3^{\lceil x_{h}\rceil}-1}{n+1}\right)
\]
as the first interval of $y_{n+1}(x)$ and so on. The infinite recursion of (\ref{eq:yn}) does not lead to an infinite loop, because we can always get the first interval of $y_n(x)$ without recursion into $y_{n+1}(x)$.

There is one further point to consider; the first interval of $M_{(n,x)}(y_{n+1}(x))$ might not be a subset of (\ref{eq:ynInterval}). Since both intervals contain $y_n(x)$, we can take their intersection as the second interval of $y_n(x)$. We might need to continue taking such intersections as long as $M_{(n,x)}(y_{n+1}(x))$ produces ambiguous terms. 

If the interval (\ref{eq:ynInterval}) is unambiguous, giving $a$ as the first term of $y_n(x)$, we can immediately produce $a$, followed by applying the arithmetic algorithm to
\[
\frac{1}{M_{(n,x)}(y_{n+1}(x)) - a} \ .
\]

\subsubsection{Validity of the Algorithm}\label{sec:validity} At this point we have given an algorithm which produces a sequence of intervals for each $y_n(x)$; it remains to prove that the sequence of intervals representing $y_1(x)$ is in fact valid as defined in section \ref{sec:FailToConverge}.
The properties of the arithmetic algorithm only guarantee that, if $x$ and $y_n(x)$ are valid sequences of intervals, then $y_{n-1}(x)$ is valid.

The algorithm already guarantees that successive ambiguous intervals will be nested; we must prove that the upper and lower bounds converge to the same point. In other words we must prove that, for any $\varepsilon > 0$, we eventually obtain a prefix $[I_0, I_1,\cdots I_k]$ of  $y_1(x)$ which confines  $y_1(x)$ within an interval of width $\leq \varepsilon$; we say that $y_1(x)$ is $\varepsilon$-\emph{bounded}.

So suppose that at some point we have ingested enough of  $y_{n+1}(x)$ to have $y_{n+1}(x)$ be $\varepsilon'$-bounded, with $y_\ell \leq y_{n+1}(x) \leq y_\ell + \varepsilon'$. The definition of $M_{(n,x)}$, along with $x_{\ell} \leq x \leq x_{h}$,  implies that  $y_n(x)$ is $\varepsilon$-bounded if
\[
\varepsilon' \leq \frac{n}{x_h}\varepsilon - \frac{x_h - x_\ell}{x_h}y_\ell
\]
But $\frac{x_h - x_\ell}{x_h}$ approaches zero as we ingest more of $x$, while $y_\ell$ is bounded above by a constant, so eventually it will suffice to have
\[
\varepsilon' \leq \frac{n}{2x_h}\varepsilon
\]

Thus $y_1(x)$ is $\varepsilon$-bounded if some $y_n(x)$ is $\varepsilon'$-bounded with
\[
\varepsilon' \leq \frac{(n-1)!}{(2x_h)^{n-1}}\varepsilon \ .
\]
Now $\frac{(n-1)!}{(2x_h)^{n-1}}$ approaches infinity, while our algorithm immediately bounds each $y_n(x)$ within an interval of width
\[
 \frac{3^{\lceil x_{h}\rceil}-1 - x_{\ell}}{n}
\]
which approaches zero; therefore some $y_n(x)$ will be $\varepsilon'$-bounded as required. The CF will converge to the correct value $e^x$, since all intervals bounding $y_1(x)$ by construction contain $e^x$.

\subsection{Logarithms} We can apply the same approach to $\log(x)$, using the series
\[
\log (x) =  2 \left( \left(\frac{x-1}{x+1}\right) + \frac{1}{3}\left(\frac{x-1}{x+1}\right)^3 + \frac{1}{5}\left(\frac{x-1}{x+1}\right)^5 \cdots + \right)
\]
which converges for all $x > 0$. Let $z=\frac{x-1}{x+1}$, and define
\[
g(w) = 1 + \frac{w}{3} + \frac{w^2}{5} + \frac{w^3}{7} + \cdots = 1 + \frac{w}{3}\left(1 + \frac{3w}{5}\left(1 + \frac{5w}{7}\left(1 +\cdots \right. \right. \right.\ .
\]
Then $\log(x) = 2zg(z^2)$, so we only need consider how to compute $g(w)$ with $0 \leq w \leq 1$.
The series for $g(w)$ can be obtained by iterating
\[
M_{(n,w)}(y) = 1 + \frac{2n-1}{2n+1}wy 
\]
i.e. the bihomographic expression $\bihomographic{2n-1}{0}{0}{2n+1}{0}{0}{0}{2n+1}$.
So we define
\[
g_n(w) = 1  + \sum_{k=1}^{\infty}\frac{2n-1}{2n+2k-1}w^k= \lim_{m\to\infty} M_{(n,w)}( M_{(n+1,w)}( \cdots M_{(m,w)}(1) \cdots ) ) \ .
\]
In particular $g_1(w) = g(w)$, and $g_n(w) = M_{(n,w)}(g_{n+1}(w))$. Since $w$ is nonnegative, we have
\begin{equation}\label{eq:logbound}
1 \leq g_n(w) \leq 1 + \frac{2n-1}{2n+1}\sum_{k=1}^{\infty} w^k = 1+\frac{2n-1}{2n+1}\frac{w}{1-w}\ .
\end{equation}
A finite prefix of the representation of $x$ will suffice to obtain a rational upper bound $B_x$ on $\frac{w}{1-w}$, so the first interval of $g_n(w)$ can be $(1, 1+\frac{2n-1}{2n+1}B_x)$, followed by applying the arithmetic algorithm to $M_{(n,w)}(g_{n+1}(w))$. However, if $\frac{2n-1}{2n+1}B_x \leq 1$, we can immediately produce 1 as the first term, followed by applying the
arithmetic algorithm to
\[
\frac{1}{M_{(n,w)}(g_{n+1}(w)) - 1}\ .
\]

\subsubsection{Validity} Proceeding along the same lines as section (\ref{sec:validity}), $g_n(w)$ will be $\varepsilon$-bounded if
$g_{n+1}(w)$ is $\frac{2n+1}{(2n-1)\hat{w}}\varepsilon$-bounded for some $w < \hat{w} < 1$; so $g_1(w)$ will be $\varepsilon$-bounded if some $g_n(w)$ is
$\frac{2n+3}{\hat{w}^{n-1}}\varepsilon$-bounded. This width approaches infinity,
while the algorithm immediately confines each $g_n(w)$ in an interval no wider than the constant $B_x$.

\subsection{Trigonometric Functions}
The cosine
\[
\cos(x) = 1 - \frac{x^2}{2!} + \frac{x^4}{4!} - \cdots = 1 - \frac{x^2}{2} \left(1 - \frac{x^2}{12}\left( 1 - \frac{x^2}{30}\left(1 + \cdots \right. \right.  \right.  
\]
can be represented by letting $w=x^2$, defining
\[
M_{n,w}(y) = 1 - \frac{wy}{2n(2n-1)}
\]
and
\begin{equation}\label{eq:cosIter}
c_n(w) = M_{n,w}(M_{n+1,w}(\cdots)) = 1 - \frac{w}{2n(2n-1)} + \frac{w^2}{(2n+2)(2n+1)2n(2n-1)} - \cdots \ .
\end{equation}
so $c_n(w) = M_{n,w}(c_{n+1}(w))$, and $c_1(w) = \cos(x)$.

If $|x| \leq \frac{\pi}{2}$, then $w < \frac{5}{2}$, and (\ref{eq:cosIter}) is an alternating series bounded between $1 - \frac{5}{4n(2n-1)}$ and 1; the CF expansion of $c_1(w)$ will start with this interval,  followed by the result of applying the arithmetic algorithm to the bihomographic expression $M_{n,w}(c_2(w))$. For $n \geq 2$, we can immediately produce 0, and we have
\[
1  \leq \frac{1}{M_{n,w}(c_{n+1}(w))} \leq \frac{4n(2n-1)}{4n(2n-1)-5} \ .
\]
This interval in turn lets us immediately produce 1; then
\[
 \frac{4n(2n-1)-5}{5} \leq \left(\frac{1}{M_{n,w}(c_{n+1}(w))} - 1 \right)^{-1} < \infty
\]
and subsequent intervals come from applying the arithmetic algorithm to
\[
\left(\frac{1}{M_{n,w}(c_{n+1}(w))} - 1\right) ^{-1} \ .
\]

For $\frac{\pi}{2} \leq |x| \leq \pi$ we use $\cos(x) = -\cos(\pi - x)$, and for all other values $\cos(x) = \cos(x \pm 2\pi)$; see section \ref{sec:spigot} for computing the CF expansion of $\pi$. We then define $\sin(x) = \cos(x-\pi/2)$ and $\tan(x)=\sin(x)/\cos(x)$.

Proof of validity is similar to the proof for $e^x$.

\subsection{Inverse Trigonometric Functions}
Finally we consider
\[
\frac{\arcsin(x)}{x} = \sum_{n=0}^{\infty} \frac{(2n)!}{4^n(n!)^2(2n+1)}x^{2n} = 1 + \frac{x^2}{6} + \frac{3x^4}{40} + \cdots
\]
with $| x | \leq 1$. Let $t_n = \frac{(2n)!}{4^n(n!)^2(2n+1)}$ and
\[
M_{n,w}(y) = 1 + \frac{t_n}{t_{n-1}}wy = 1 + \frac{(2n-1)^2}{2n(2n+1)}wy \ .
\]
Define $w=x^2$ and
\[
a_n(w) = M_{n,w}(M_{n+1,w}(\cdots)) = 1 + \frac{t_n}{t_{n-1}}w+ \frac{t_{n+1}}{t_{n-1}}w^2 + \cdots
\]
so $\arcsin(x) = x a_1(w)$. Since $0 \leq w \leq 1$, and $\{t_n\}$ is strictly decreasing, we can obtain the initial interval of $a_n(w)$ from
\[
1 \leq a_n(w) \leq 1 + \frac{(\arcsin(1)-1)w}{t_{n-1}} < 1 + \frac{4w}{7t_{n-1}} \ .
\]
For $n=1$ we have $t_0=1$ and can immediately produce 1 as the first term.
Proof of validity is similar to the proof for logarithm.

\subsection{Spigot Algorithms}\label{sec:spigot}
Expressions like (\ref{eq:spigot}), representing an infinite sum as an infinitely nested sequence of rational functions, are central to the spigot algorithm\cite{Gibbons2016,Rabinowitz2016} for computing digits of $\pi$. One version of that algorithm uses Gosper's ``accelerated" series \cite{acceleration}:
\[
\pi = 3 + \frac{1 \times 1}{3 \times 4 \times 5}\big(8 + \frac{2 \times 3}{3 \times 7 \times 8}\big( \cdots 5i-2 + \frac{i(2i-1)}{3(3i+1)(3i+2)} \cdots \big) \big) \ .
\]
We can define
\[
M_i(x) = \cdots 5i-2 + \frac{i(2i-1)}{3(3i+1)(3i+2)}x
\]
and $p_i = \lim_{i \to \infty} M_i(M_{i+1}(\cdots))$;
then $\pi = p_1$. The spigot algorithm uses the fact that
\begin{equation}\label{eq:pi}
\frac{27i-12}{5} \leq p_i \leq \frac{27}{5}i - \frac{216}{125} \ .
\end{equation}
Following the same pattern as used for transcendental functions, we could obtain the CF terms of $\pi = p_1$ using $p_i = M_n(p_{i+1})$, taking (\ref{eq:pi}) as the first interval of $p_i$. The spigot algorithm uses the same equation to produce digits; when the floor of the current rational expression is known to be $k$, we transform it as $M \leftarrow 10(M-k)$ to obtain an expression whose floor is the next digit.

In fact, in the same way, we could produce base-$b$ digits of $e^x$ (with $x$ given as a CF!) by redefining $\mbox{produce}(a, M)$ as
\[
M(x,y) \leftarrow b(M(x,y) - a)
\]
when we compute $y_1(x)$;
after such a transformation, $M(x,y)$ represents a real number whose floor is the next digit, rather than the next CF term. 
For that matter, if $x$ were given as a sequence of digits, we could read it by redefining $\mbox{ingest\_x}(a,M)$ as the transformation
\[
x \leftarrow a + \frac{x}{b} \ .
\]
 After this transformation, $M(x,y)$ represents the same value, with $x=[a,x_1,x_2\cdots]$ replaced by $x=[x_1,x_2\cdots]$; the upper and lower bounds on $M(x,y)$ would be determined by $0 \leq x < b$ instead of $1 <  x \leq \infty$.

\subsection{Related Work}
Davison \cite{DavisonCFexp} gives an algorithm for the CF terms of $e^{\alpha}$ with $\alpha$ rational;
it computes convergents using the fact that if
\[
\begin{pmatrix}
p_n & q_n \\
r_n & s_n \end{pmatrix} = \prod_{k=0}^n \begin{pmatrix} (2k+1)m+\ell & (2k+1)m \\
                                                       (2k+1)m &  (2k+1)m-\ell \end{pmatrix}
\]
then $p_n/r_n$ and $q_n/s_n$ approach $e^{\ell/m}$.

We note that standard computational software packages such as Mathematica\footnote{\href{wolfram.com}{wolfram.com}} and Pari\footnote{\href{https://pari.math.u-bordeaux.fr}{https://pari.math.u-bordeaux.fr}}\footnote{\href{https://wims.univ-cotedazur.fr/~wims/wims.cgi}{https://wims.univ-cotedazur.fr/~wims/wims.cgi}} can compute CF expansions of complex expressions to high accuracy very much faster than our Haskell code.
Our algorithms, however, can operate on real numbers defined by arbitrary sequences of CF terms $[a_0,a_1,\cdots]$; more generally, on any number defined by converging sequences of rational upper and lower bounds.
The algorithms could also be adapted to other functions defined by Taylor series, if we can obtain appropriate bounds in the same manner as (\ref{eq:ynInterval}) and (\ref{eq:logbound}), along with corresponding proofs of convergence.

\section{Optimizing Computation of $e^x$ and $\log(x)$}
In section (\ref{sec:validity}), the larger $x$ is, the larger $n$ will have to be before $y_n(x)$ is bounded in  a given interval. In the algorithm, this corresponds to recursing more deeply into $M_{(1,x)}( M_{(2,x)}( \cdots ) )$ before we get tight bounds on $y_1(x)$. Thus our Haskell implementation computes $y_1(x)$ only for $0 \leq x < 1$; when $x>1$ we use
\[
e^x = e^{\left\lfloor x \right\rfloor}e^{x - \left\lfloor x \right\rfloor} \ .
\]
Ordinary multiplication can be used to obtain $e^{\left\lfloor x \right\rfloor}$, since
\[
e = [2,1,2,1,1,4,1,1,6,1,1,8,1,1,10,\cdots] \ .
\]
Indeed, when the fractional part of $x$ is greater than $\frac{1}{2}$ we could use
\[
\sqrt{e} = [1,1,1,1,5,1,1,9,1,1,13\cdots]
\]
and $e^x = e^{\left\lfloor x \right\rfloor }e^{x - \left\lfloor x \right\rfloor- \frac{1}{2}}\sqrt{e}$.

Similarly, computation of $g(w)$ recurses less deeply when $w$ is small, i.e. when $x$ is close to 1; for $x>e$ our Haskell implementation uses $\log(x) = 1 + \log(x/e)$.

\bibliographystyle{plain}
\bibliography{cfAlgorithm}

\end{document}